\documentclass{amsart}

\usepackage{amsmath,amssymb,amscd,amsfonts}

\newtheorem{theorem}{Theorem}
\newcommand{\bt}{\begin{theorem}}
\newcommand{\et}{\end{theorem}}
\newtheorem{lemma}{Lemma}
\newcommand{\bl}{\begin{lemma}}
\newcommand{\el}{\end{lemma}}
\newtheorem{corollary}{Corollary}
\newcommand{\bc}{\begin{corollary}}
\newcommand{\ec}{\end{corollary}}
\newcommand{\beq}{\begin{equation}}
\newcommand{\eeq}{\end{equation}}
\newcommand{\benum}{\begin{enumerate}}
\newcommand{\eenum}{\end{enumerate}}

\newcommand{\R}{\ensuremath{\mathbf R}}

\newcommand{\mcb}{\ensuremath{ \mathcal B}}

\newcommand{\mce}{\ensuremath{ \mathcal E}}
\newcommand{\mcf}{\ensuremath{ \mathcal F}}

\newcommand{\mba}{\ensuremath{ \mathbf a}}
\newcommand{\mbb}{\ensuremath{ \mathbf b}}
\newcommand{\mbc}{\ensuremath{ \mathbf c}}
\newcommand{\mbe}{\ensuremath{ \mathbf e}}

\newcommand{\mbf}{\ensuremath{ \mathbf f}}
\newcommand{\mbu}{\ensuremath{ \mathbf u}}
\newcommand{\mbv}{\ensuremath{ \mathbf v}}

\newcommand{\mbx}{\ensuremath{ \mathbf x}}
\newcommand{\mby}{\ensuremath{ \mathbf y}}

\newcommand{\Rn}{\ensuremath{ \mathbf{R}^n }}
\DeclareMathOperator{\colsum}{\text{colsum}}
\DeclareMathOperator{\rowsum}{\text{rowsum}}
\newcommand{\bsmallmat}{\left(\begin{smallmatrix}}
\newcommand{\esmallmat}{\end{smallmatrix}\right)}
\newcommand{\mbj}{\ensuremath{\mathbf j}}

\DeclareMathOperator{\card}{\text{card}}
\DeclareMathOperator{\conv}{\text{conv}}

\newcommand{\bmat}{\left(\begin{matrix}}
\newcommand{\emat}{\end{matrix}\right)}
\DeclareMathOperator{\qand}{\quad\text{and}\quad}
\DeclareMathOperator{\qqand}{\qquad\text{and}\qquad}

\DeclareMathOperator{\vectoran}{\left( \begin{matrix} a_1 \\ \vdots \\ a_n \end{matrix}\right)}

\DeclareMathOperator{\vectorxn}{\left( \begin{matrix} x_1 \\ \vdots \\ x_n \end{matrix}\right)}
\DeclareMathOperator{\vectoryn}{\left( \begin{matrix} y_1 \\ \vdots \\ y_n \end{matrix}\right)}

\DeclareMathOperator{\vectorsmallan}{\left( \begin{smallmatrix} a_1 \\ \vdots \\ a_n \end{smallmatrix}\right)}
\DeclareMathOperator{\vectorsmallbn}{\left( \begin{smallmatrix} b_1 \\ \vdots \\ b_n \end{smallmatrix}\right)}
\DeclareMathOperator{\vectorsmallcn}{\left( \begin{smallmatrix} c_1 \\ \vdots \\ c_n \end{smallmatrix}\right)}

\DeclareMathOperator{\vectorsmallun}{\left( \begin{smallmatrix} u_1 \\ \vdots \\ u_n \end{smallmatrix}\right)}
\DeclareMathOperator{\vectorsmallvn}{\left( \begin{smallmatrix} v_1 \\ \vdots \\ v_n \end{smallmatrix}\right)}

\DeclareMathOperator{\vectorsmallxn}{\left( \begin{smallmatrix} x_1 \\ \vdots \\ x_n \end{smallmatrix}\right)}

\DeclareMathOperator{\vectorsmallx2}{\left( \begin{smallmatrix} x_1 \\ x_2 \end{smallmatrix}\right)} 
\DeclareMathOperator{\vectorsmallb2}{\left( \begin{smallmatrix} b_1 \\ b_2 \end{smallmatrix}\right)}

\DeclareMathOperator{\col}{\text{col}}

\title{The Muirhead-Rado inequality, 2:  Symmetric means and inequalities}
\author{Melvyn B. Nathanson}
\address{Lehman College (CUNY), Bronx, New York 10468} 
\email{melvyn.nathanson@lehman.cuny.edu}



\thanks{Supported in part by PSC CUNY Grant \# 66197-00 54.}

\begin{document}

\begin{abstract}
This paper gives elementary proofs of  the Muirhead and Rado inequalities.
\end{abstract}

\maketitle

\section{Symmetric means of monomial functions}

A \emph{nonnegative vector}\index{vector!nonnegative} is a vector with nonnegative coordinates.
The \emph{nonnegative octant}\index{octant!nonnegative} in \Rn\ is 
\[
\Rn_{\geq 0} = 
\left\{ \mba = \vectoran \in \Rn : \text{$a_i \geq 0$ for all $i \in \{1,\ldots, n\}$} \right\}.
\]
A \emph{positive vector}\index{vector!positive} is a vector with positive coordinates.
The \emph{positive octant}\index{vector!positive}\index{octant!positive} in \Rn\ is 
\[
\Rn_{>0} = 
\left\{ \mbx = \vectorxn \in \Rn : \text{$x_i > 0$ for all $i \in \{1,\ldots, n\}$} \right\}.
\]  

Every nonnegative vector $\mba$ defines the  \emph{monomial function}\index{monomial function}   
$\mbx^{\mba}$ on the positive octant as follows:  
For $\mba = \vectorsmallan \in \R^n_{\geq 0}$ and  $\mbx = \vectorsmallxn \in \R^n_{> 0}$,     
let
\[
\mbx^{\mba} = x_1^{a_1} x_2^{a_2} \cdots x_n^{a_n}.
\]  
We call \mba\   the \emph{exponent vector}\index{exponent vector} 
of the monomial $\mbx^{\mba}$.  Note that $\mbx^{\mba} > 0$ for all $\mbx \in \R^n_{>0}$. 

Let $S_n$ be the symmetric group.  
The \emph{symmetric mean} of the monomial  function $\mbx^{\mba} $ is the function 
\[
[\mbx^{\mba} ]_{S_n} =  \frac{1}{n!} \sum_{\sigma \in S_n}  
x_{\sigma(1)}^{a_1}  x_{\sigma(2)}^{a_2}  \cdots x_{\sigma(n)}^{a_n}.
\]
For every subgroup $G$ of $S_n$, the \emph{$G$-symmetric mean} 
of the monomial $\mbx^{\mba} $ is the function 
\[
[\mbx^{\mba} ]_{G} =  \frac{1}{|G|} \sum_{\sigma \in G}  
x_{\sigma(1)}^{a_1}  x_{\sigma(2)}^{a_2}  \cdots x_{\sigma(n)}^{a_n}.
\]

Let $\sigma \in S_n$.  If  $i \in \{1,\ldots, n\}$ and $j = \sigma(i)$, then  
\[
x_{\sigma(i)}^{a_i} =  x_j^{ a_{ \sigma^{-1}(j) }}  
\]
and  
\[
\prod_{i=1}^n x_{\sigma(i)}^{a_i} = \prod_{j=1}^n x_j^{ a_{\sigma^{-1}(j)}}. 
\]
Because $\{\sigma \in G \} = \{ \sigma^{-1} \in G \}$, we have 
\begin{align*} 
[\mbx^{\mba} ]_{G} 
& =  \frac{1}{|G|} \sum_{\sigma \in G}  \prod_{i=1}^n x_{\sigma(i)}^{a_i} 
 =  \frac{1}{|G|} \sum_{\sigma \in G}  \prod_{j=1}^n x_j^{ a_{\sigma^{-1}(j)}} \\
& =  \frac{1}{|G|} \sum_{\sigma \in G}   \prod_{j=1}^n x_j^{ a_{\sigma(j)}} 
 =  \frac{1}{|G|} \sum_{\sigma \in G}   
x_1^{ a_{ \sigma(1) }}  x_2^{ a_{ \sigma(2) }}  \cdots x_n^{ a_{ \sigma(n) }} .
\end{align*}
If $\mbx  \in \R^n_{> 0}$ is the constant vector with all coordinates equal to $x$, then 
\beq              \label{Muirhead:constant}
[\mbx^{\mba}]_G 
=  \frac{1}{|G|} \sum_{\sigma \in G}   \prod_{j=1}^n x^{ a_{\sigma(j)}} 
=  \frac{1}{|G|} \sum_{\sigma \in G}  x^{\sum_{i=1}^n a_i} 
= x^{\sum_{i=1}^n a_i}.
\eeq

Let $G$ be a subgroup of $S_n$ and let $\mba = \vectorsmallan$ and $\mbb = \vectorsmallbn$ 
be nonnegative vectors such that 
\[
\sum_{j=1}^n a_j = \sum_{j=1}^n b_j.
\]
If 
\beq        \label{Muirhead:MonomialIneq}
\left[\mbx^{\mbb}\right]_G < \left[\mbx^{\mba}\right]_G 
\eeq
for all vectors $\mbx \in \R^n_{>0}$, 
then~\eqref{Muirhead:MonomialIneq} is called a \emph{monomial inequality 
with respect to the subgroup $G$} determined by \mba\ and \mbb. 
A  \emph{monomial inequality} is a   monomial inequality 
with respect to the symmetric group $S_n$. 
Here are some examples.  
If $\mba = \bsmallmat 7 \\ 3 \esmallmat  \in \R^2_{\geq 0}$, 
then $\mbx^{\mba} = x_1^7 x_2^3$.  
The $S_2$-symmetric mean of $\mbx^{\mba} $ is 
\[
[\mbx^{\mba} ]_{S_2} =  \frac{1}{2}\left( x_1^7 x_2^3 + x_2^7 x_1^3  \right) 
=  \frac{1}{2}\left( x_1^7 x_2^3 + x_1^3 x_2^7 \right). 
\]
 If $\mbb = \bsmallmat 6 \\ 4 \esmallmat  \in \R^2_{\geq 0}$, 
then $\mbx^{\mbb} = x_1^6 x_2^4$ and 
the $S_2$-symmetric mean of $\mbx^{\mbb} $ is 
\[
[\mbx^{\mbb} ]_{S_2} =  \frac{1}{2}\left( x_1^6 x_2^4 + x_1^4 x_2^6 \right). 
\] 
It is straightforward to check that 
\beq            \label{Muirhead:73-64} 
x_1^6 x_2^4 + x_1^4 x_2^6 <  x_1^7 x_2^3 + x_1^3 x_2^7  
\eeq 
for all nonconstant positive vectors 
$\mbx = \left( \begin{smallmatrix} x_1 \\ x_2 \end{smallmatrix}\right)$  
and so $[\mbx^{\mbb} ]_{S_2}  < [\mbx^{\mba} ]_{S_2} $ is a monomial inequality.  

Another example. 
Associated to the vector $\mbe_1 = \bsmallmat 1\\ 0 \\ \vdots \\ 0 \esmallmat \in \Rn_{\geq 0}$ 
is the monomial function 
\[
\mbx^{\mbe_1} = x_1^1  x_2^0  \cdots x_n^0= x_1.
\]
For every ordered pair of integers $(i,j)$ with $i,j \in \{1,\ldots, n\}$, there are $(n-1)!$ permutations 
$\sigma \in S_n$ with $\sigma(i) = j$.  
In particular, for all $j \in \{1,\ldots, n\}$, there are $(n-1)!$ permutations 
$\sigma \in S_n$ with $\sigma(1) = j$.  For all $\mbx = \vectorsmallxn \in \R^n_{>0}$, we have  
\begin{align*}
[\mbx^{\mbe_1} ]_{S_n}  & =  \frac{1}{n!} \sum_{\sigma \in S_n}  x_{\sigma(1)}^1  x_{\sigma(2)}^0  \cdots x_{\sigma(n)}^0 \\
& =  \frac{1}{n!} \sum_{\sigma \in S_n}  x_{\sigma(1)} 
=  \frac{1}{n!} \sum_{j =1}^n  (n-1)!x_j \\
& =  \frac{1}{n} \sum_{j =1}^n x_j.   
\end{align*}
Thus, the symmetric mean of the monomial $\mbx^{\mbe_1}$ 
is the \index{arithmetic mean}\emph{arithmetic mean} of the 
positive numbers $x_1,\ldots, x_n$.  

Associated to the vector $\mbf_n  = \bsmallmat 1/n \\ 1/n \\ \vdots \\ 1/n \esmallmat \in \Rn_{\geq 0}$ 
is the monomial function 
\[
\mbx^{\mbf_n} =  x_1^{1/n}  x_2^{1/n}    \cdots x_n^{1/n} 
=  \left(x_1 x_2 \cdots x_n \right)^{1/n}  
\]
which is is the \index{geometric mean}\emph{geometric mean} of the 
positive numbers $x_1,\ldots, x_n$.  Because 
$x_{\sigma(1)} x_{\sigma(2)}   \cdots x_{\sigma(n)} = x_1 x_2 \cdots x_n$ 
for all $\sigma \in S_n$,   the symmetric mean 
\begin{align*}
[\mbx^{\mbf_n}]_{S_n} 
& =   \frac{1}{n!} \sum_{\sigma \in S_n}  \left(x_{\sigma(1)} x_{\sigma(2)}   \cdots x_{\sigma(n)}\right)^{1/n}   \\
& =   \frac{1}{n!} \sum_{\sigma \in S_n}  \left(x_1 x_2 \cdots x_n \right)^{1/n}   \\
& =    \left(x_1 x_2 \cdots x_n \right)^{1/n}  
\end{align*}
is also the geometric mean of the numbers $x_1,\ldots, x_n$.  
The arithmetic and geometric mean inequality is equivalent to the  monomial inequality
\[
[\mbx^{\mbf_n}]_{S_n} < [\mbx^{\mbe_1}]_{S_n}  
\]
for all nonconstant vectors $\vectorsmallxn \in \Rn_{>0}$.  

There are infinitely many ordered pairs $( \mba, \mbb)$ of nonnegative vectors such that 
\[
[\mbx^{\mbb}]_{S_n} < [\mbx^{\mba}]_{S_n} 
\]
 for all nonconstant vectors $\mbx  \in \Rn_{>0}$.  
 Muirhead's inequality (Theorem~\ref{Muirhead:theorem:Muirhead}) 
and its converse (Theorem~\ref{Muirhead:theorem:Muirhead-2}) 
classify the pairs of vectors \mba\ and \mbb\ 
that determine monomial inequalities.   Richard Rado 
(Theorems~\ref{Muirhead:theorem:Rado-1} and~\ref{Muirhead:theorem:Rado-2}) 
classified monomial inequalities with respect to a subgroup $G$ of $S_n$.  
The purpose of this paper is to prove these results.

\section{Simple inequalities}
Here is a trivial fact, from which we obtain infinitely many monomial inequalities. 

\bt         \label{Muirhead:theorem:SimpleIneq}
The product of two   real numbers is positive if and only if both numbers are positive 
or  both numbers are negative. 
\et

\bc                    \label{Muirhead:corollary:SimpleIneq-1} 
Let $\mba =  \left( \begin{smallmatrix} a_1 \\ a_2 \end{smallmatrix}  \right) \in \R^2_{> 0}$ and 
$\mbx = \left( \begin{smallmatrix} x_1 \\  x_2  \end{smallmatrix} \right) \in \R^2_{>0}$.   
If $x_1 \neq x_2$, then 
\beq                    \label{Muirhead:SimpleIneq-ab}
x_1^{a_1} x_2^{a_2} + x_1^{a_2}x_2^{a_1}  < x_1^{a_1+a_2} + x_2^{a_1+a_2}. 
\eeq
If $x_1 = x_2$, then 
\[
x_1^{a_1}x_2^{a_2} + x_1^{a_2}x_2^{a_1}  = x_1^{a_1+a_2} + x_2^{a_1+a_2}. 
\]
\ec

\begin{proof}
If $0 < x_2 < x_1$, then $x_1^{a_1}-x_2^{a_1} > 0$ and $x_1^{a_2}-x_2^{a_2} > 0$.   
If $0 < x_1 < x_2$, then $x_1^{a_1}-x_2^{a_1} < 0$ and $x_1^{a_2}-x_2^{a_2} < 0$.  
Thus, the numbers $x_1^{a_1}-x_2^{a_1}$ and $x_1^{a_2}-x_2^{a_2}$ 
are both positive or both negative.  
Applying Theorem~\ref{Muirhead:theorem:SimpleIneq} 
to these numbers gives 
\[
x_1^{a_1+a_2} + x_2^{a_1+a_2} - x_1^{a_1}x_2^{a_2} - x_1^{a_2}x_2^{a_1} = (x_1^{a_1}-x_2^{a_1})(x_1^{a_2}-x_2^{a_2}) > 0.  
\]
This proves~\eqref{Muirhead:SimpleIneq-ab}.  If $x_1 = x_2$, then 
\[
x_1^{a_1}x_2^{a_2} + x_1^bx_2^{a_1}  = x_1^{a_1+a_2} + x_2^{a_1+a_2}. 
\]
This completes the proof. 
\end{proof}

\bc                  \label{Muirhead:corollary:SimpleIneq-2}
Let $\rho$, $ \delta $, and $ \Delta$ satisfy the inequalities 
\[
\rho > 0 \qqand 
| \delta | < \Delta. 
\]
Let $\mbx = \left( \begin{smallmatrix} x_1 \\  x_2  \end{smallmatrix}\right) \in \R^2_{>0}$. 
If $x_1 \neq x_2$, then 
\beq                    \label{Muirhead:SimpleIneq-3}
x_1^{ \rho + \delta} x_2^{  \rho  - \delta} + x_1^{  \rho  -\delta} x_2^{ \rho  + \delta} 
< x_1^{ \rho + \Delta} x_2^{ \rho   -\Delta} + x_1^{ \rho   -\Delta} x_2^{  \rho +\Delta}.
\eeq
\ec

\begin{proof} 
We have $\Delta + \delta > 0$ and $\Delta - \delta > 0$.
Applying Corollary~\ref{Muirhead:corollary:SimpleIneq-1} 
with $\mba = \bmat \Delta + \delta \\ \Delta - \delta \emat   \in \R^2_{>0} $ gives 
\[
x_1^{\Delta + \delta}  x_2^{\Delta - \delta}  + x_1^{\Delta - \delta} x_2^{\Delta + \delta}   
< x_1^{2\Delta} + x_2^{2\Delta }.  
\]
Multiplying this inequality by $(x_1x_2)^{-\Delta}$ gives 
\[
x_1^{\delta} x_2^{-\delta} + x_1^{ -\delta} x_2^{ \delta} 
< x_1^{ \Delta} x_2^{ -\Delta} + x_1^{ -\Delta} x_2^{ \Delta}.  
\]  
Multiplying  by $(x_1x_2)^{\rho}$ completes the proof. 
\end{proof}

For example, letting $\rho = 5$, $\delta = 1$,  and $\Delta = 2$
in inequality~\eqref{Muirhead:SimpleIneq-3} gives  inequality~\eqref{Muirhead:73-64}. 
Similarly,  for all nonconstant vectors $\mbx \in \R^2_{>0}$,  we have 
the chain of inequalities 
\begin{align*}
2x_1^5 x_2^5 
& < x_1^6 x_2^4 +  x_1^4 x_2^6 < x_1^7 x_2^3 + x_1^3 x_2^7 \\ 
& < x_1^8 x_2^2 +  x_1^2 x_2^8 < x_1^9 x_2 + x_1 x_2^9 
 < x_1^{10} +  x_2^{10}.
\end{align*}
These inequalities fail if $x_1=1$ and $x_2 = -1$ 
or, more generally, if $x_1x_2 \leq 0$.

\bt                    \label{Muirhead:theorem:Muirhead-0}
Let vectors $\mba = \left( \begin{smallmatrix} a_1 \\ a_2 \end{smallmatrix} \right) \in \R^2_{\geq 0}$ 
and $\mbb = \vectorsmallb2  \in \R^2_{\geq 0}$ satisfy 
\[
 a_2 < b_2 \leq b_1 < a_1 \qqand a_1+a_2 =  b_1+b_2. 
\]
If  $\mbx = \vectorsmallx2  \in \R^2_{> 0}$  with $x_1 \neq x_2$, then 
$[\mbx]_{S_2}^{\mbb} < [\mbx]_{S_2}^{\mba}$ or, equivalently, 
\beq                                       \label{Muirhead:Muirhead-0}
x_1^{ b_1} x_2^{ b_2} + x_1^{ b_2} x_2^{b_1} 
< x_1^{ a_1} x_2^{ a_2} + x_1^{a_2} x_2^{ a_1}.  
\eeq   
\et

\begin{proof}
Let 
\[
\rho = \frac{a_1+a_2}{2} = \frac{b_1+b_2}{2} 
\] 
and 
\[
\delta = \frac{b_1-b_2}{2} \qqand \Delta = \frac{a_1-a_2}{2}. 
\]
We have   
\[
\rho > 0 \qqand  0 \leq \delta < \Delta.  
\]
Applying inequality~\eqref{Muirhead:SimpleIneq-3} with 
\[
\rho+\delta = b_1 \qqand \rho - \delta = b_2
\]
and 
\[
\rho+\Delta = a_1 \qqand \rho - \Delta = a_2
\]
gives inequality~\eqref{Muirhead:Muirhead-0}.    
This completes the proof. 
\end{proof}

This is Muirhead's inequality for monomials in two variables.  
Inequality~\eqref{Muirhead:73-64} is the special case 
$\mba= \bsmallmat  7 \\ 3 \esmallmat$ and   
$\mbb = \bsmallmat  6 \\ 4 \esmallmat$.

\section{Symmetric means of functions of $n$ variables} 

Recall the \emph{Kronecker delta} 
\[
\delta_{i,j} 
= \begin{cases}
1 & \text{if $i = j$  } \\
0 & \text{if $i \neq j$ .}
\end{cases}
\]
The \emph{standard basis} for \Rn\ is   $\mce = \{\mbe_1,\ldots, \mbe_n\}$, 
where $\mbe_j = \bsmallmat \delta_{1,j} \\ \vdots \\ \delta_{n,j} \esmallmat$. 
Every permutation $\sigma \in S_n$ defines a linear operator on \Rn\ by 
$\sigma ( \mbe_j) = \mbe_{\sigma(j)}$.  
If $\mbx =  \sum_{j=1}^n x_j \mbe_j \in \R^n$, then 
\[
\sigma (\mbx ) 
=   \sum_{j=1}^n x_j \sigma \left(\mbe_j  \right) 
=   \sum_{j=1}^n x_j  \mbe_{\sigma(j)}
=   \sum_{j=1}^n x_{\sigma^{-1}(j)}  \mbe_j 
\]
and so 
\beq                \label{Muirhead: sigma-z}
\sigma \vectorxn
= \bmat x_{\sigma^{-1}(1)} \\ \vdots \\ x_{\sigma^{-1}(n)} \emat.
\eeq

Let $\Omega$ be a subset of \Rn\ that is closed under the action of $S_n$, 
that is, $\mbx \in \Omega$ implies $\sigma ( \mbx ) \in \Omega$ for all $\sigma \in S_n$. 
Let $\mcf(\Omega)$ be the set of real-valued functions defined on $\Omega$.  

For every function $f \in \mcf(\Omega)$ and every permutation $\sigma$ 
in the symmetric group $S_n$, define the function $\sigma f \in \mcf(\Omega)$ by 
\beq                                 \label{Muirhead:actionF(Omega)}
(\sigma f) (\mbx) = f\left(\sigma^{-1}( \mbx) \right).  
\eeq
For all $\sigma,\tau \in S_n$ we have 
\begin{align*}
 \left( \tau (\sigma f)  \right) (\mbx) & = (\sigma f)\left( \tau^{-1}( \mbx) \right) 
=  f \left( \sigma^{-1}\left( \tau^{-1} ( \mbx )  \right) \right) \\
& =  f \left( \left(\sigma^{-1} \tau^{-1} \right)( \mbx )   \right) =  f \left( \left( \tau\sigma \right)^{-1} ( \mbx )   \right) \\
& =   \left( \left( \tau\sigma \right) f  \right) \left(  \mbx  \right). 
\end{align*} 
Thus, 
\[
 \left( \tau\sigma \right) f =  \tau (\sigma f) 
\]
and~\eqref{Muirhead:actionF(Omega)} defines an action of the group $S_n$ 
on the function space $\mcf(\Omega)$.

Let $G$ be a subgroup of $S_n$ of order $|G|$.  
The \emph{$G$-symmetric mean}\index{symmetric mean}\index{symmetric mean!$G$-symmetric mean} 
of the function $f$ 
is the   function $[f]_G \in \mcf(\Omega)$ defined by 
\[
[f]_G = \frac{1}{|G|} \sum_{\sigma \in G} \sigma f.
\]
Because $G = \{ \sigma^{-1}: \sigma \in G\}$, we have 
 \begin{align*}
[f]_G(\mbx) & = \frac{1}{|G|} \sum_{\sigma \in G} (\sigma f)(\mbx) 
=  \frac{1}{|G|} \sum_{\sigma \in G} f\left(\sigma^{-1}( \mbx) \right) \\ 
& =  \frac{1}{|G|} \sum_{\sigma \in G} f\left(\sigma( \mbx) \right).
 \end{align*}
The \emph{symmetric mean} of the function $f$ 
is the $S_n$-symmetric mean.

\bl                             \label{Muirhead:lemma:symmetrization}
Let $G$ be a subgroup of $S_n$ and let $\tau$ in $G$.  
For all functions $f \in \mcf(\Omega)$, 
\[
[\tau f]_{G} = [f]_{G}
\]
\el

\begin{proof}
From the identity 
\[
 G \tau  = \{\sigma \tau: \sigma \in G \} =  G
\]
we obtain  
\begin{align*}
[\tau f]_{G} & = \frac{1}{|G|} \sum_{\sigma \in G} \sigma(\tau f) 
= \frac{1}{|G|} \sum_{\sigma \in G} (\sigma\tau )f  = \frac{1}{|G|} \sum_{\sigma \in G} \sigma f  = [ f]_{G}.
\end{align*} 
This completes the proof. 
\end{proof}

\bl           \label{Muirhead:lemma:a-permute} 
Let $G$ be a subgroup of $S_n$.
For every nonnegative vector $\mba = \vectorsmallan  \R^n_{\geq 0}$ 
and every permutation $\tau \in G$, 
\[
[\mbx^{\mba} ]_G = [\mbx^{\tau( \mba ) } ]_G.
\]
\el

\begin{proof}
From~\eqref{Muirhead: sigma-z} we have  
\[
\tau ( \mba )  = \bmat  a_{\tau^{-1}(1)} \\ \vdots \\ a_{\tau^{-1}(n) } \emat 
\qqand
\mbx^{\tau ( \mba ) } = x_1^{a_{\tau^{-1} (1)}}    \cdots x_n^{a_{\tau^{-1}(n)}}. 
\]
Because $G  \tau = \{\sigma\tau:\sigma\in G\} = G$  for all $\tau \in G$, we have  
\begin{align*}
[\mbx^{\tau( \mba ) } ]_G 
& = \frac{1}{n!} \sum_{\sigma \in G} x_{\sigma(1)}^{a_{\tau^{-1} (1)}}   
 \cdots x_{\sigma(n)}^{a_{\tau^{-1}(n)}} \\
 & = \frac{1}{n!} \sum_{\sigma \in G} x_1^{a_{\tau^{-1} \sigma^{-1} (1)}}  
\cdots x_n^{a_{\tau^{-1} \sigma^{-1} (n)}}  \\ 
  & = \frac{1}{n!} \sum_{\sigma \in G} x_1^{a_{ ( \sigma\tau)^{-1} (1)}}    
  \cdots x_n^{a_{ ( \sigma\tau)^{-1} (n)}}  \\
 & = \frac{1}{n!} \sum_{\sigma \in G} x_1^{a_{ \sigma^{-1} (1)}}  
  \cdots x_n^{a_{ \sigma^{-1}  (n)}}  \\
 & = [\mbx^{\mba} ]_G .
\end{align*} 
This completes the proof. 
\end{proof}


\bl           \label{Muirhead:lemma:transposition} 
Let $G$ be a subgroup of $S_n$ that contains the transposition $\tau = (j,k)$ with $j < k$.  
For every nonnegative vector $\mba = \vectorsmallan \in \R^n_{\geq 0}$ 
and positive vector $\mbx = \vectorsmallxn \in\R^n_{>0}$, 
\[ 
[\mbx^{\mba}]_{G}  
 = \frac{1}{2|G|} \sum_{\sigma \in G} 
 \left( x_{\sigma  (j)}^{a_j} x_{\sigma (k)}^{a_{k} }  +  x_{\sigma  (j)}^{a_{k}}  x_{\sigma (k)}^{a_{j} }  \right) 
\prod_{\substack{i=1 \\ i \neq j,k }}^n  x_{\sigma(i)}^{a_i}. 
\]
\el

\begin{proof}
The transposition $\rho = (j,k)$ acts  as follows:  
\[
\text{if} \quad \mba = \bsmallmat a_1 \\ \vdots \\ a_{j-1} \\a_j \\a_{j+1} \\ \vdots \\  a_{k-1} 
\\a_{k} \\a_{k+1}\\ \vdots \\ a_n \esmallmat
\qquad\text{then} \quad 
\rho( \mba )  = \bsmallmat a_1 \\ \vdots \\ a_{j-1} \\a_{k} \\a_{j+1} \\ \vdots \\  
a_{k-1} \\a_j \\a_{k+1}\\ \vdots \\ a_n \esmallmat.
\]  
By Lemma~\ref{Muirhead:lemma:a-permute}, 
\begin{align*}
2[\mbx^{\mba}]_{G}  
& = [\mbx^{\mba}]_{G}  + [\mbx^{\tau( \mba ) }]_{G} \\
& =  \frac{1}{|G|} \sum_{\sigma \in G} x_{\sigma  (1)}^{a_1} \cdots  x_{\sigma  (j)}^{a_j} 
\cdots x_{\sigma (k)}^{a_{k} } 
 \cdots x_{\sigma (n)}^{a_n} \\
 & \hspace{1.5cm} + 
  \frac{1}{|G|} \sum_{\sigma \in G} x_{\sigma  (1)}^{a_1} \cdots  x_{\sigma  (j)}^{a_{k}} 
\cdots x_{\sigma (k)}^{a_{j} } 
 \cdots x_{\sigma (n)}^{a_n}   \\ 
 & = \frac{1}{|G|} \sum_{\sigma \in G} 
 \left( x_{\sigma  (j)}^{a_j} x_{\sigma (k)}^{a_{k} }  +  x_{\sigma  (j)}^{a_{k}}  x_{\sigma (k)}^{a_{j} }  \right) 
\prod_{\substack{i=1 \\ i \neq j,k }}^n  x_{\sigma(i)}^{a_i}. 
\end{align*}
Dividing by 2 completes the proof. 
\end{proof}

\bl                           \label{Muirhead:lemma:Muirhead-transposition}
Let $G$ be a subgroup of $S_n$ that contains the transposition $\tau = (j,k)$  with $j < k$.  
Let $\mbu = \vectorsmallun$ be a  nonnegative vector such that $u_j > u_k$ and let 
\[
\rho = \frac{u_j+u_k}{2} \qqand \Delta =  \frac{u_j - u_k}{2}.
\]
Then  
\[
0 < \Delta \leq \rho.
\] 
Let
\[
0 \leq \delta < \Delta.    
\]
Define the  nonnegative vector $\mbv =  \vectorsmallvn$ by 
\[
v_j = \rho + \delta 
\qqand 
v_k = \rho - \delta
\]
and
\[
v_i = u_i \qquad\text{if $i \neq j,k$}.
\]
If  $\mbx = \vectorsmallxn \in \Rn_{>0}$   is a nonconstant vector such that 
$x_{\sigma(j)} \neq x_{\sigma(k)}$ for some $\sigma \in G$, then 
\[
[\mbx^{\mbv}]_G < [\mbx^{\mbu}]_G.   
\]      
\el

\begin{proof}
We have $v_i = u_i \geq 0$ for $i \neq j,k$  and 
\[
0 \leq  u_k = \rho - \Delta < \rho - \delta = v_k \leq v_j =  \rho + \delta  < \rho + \Delta = u_j
\]
and so the vector is  \mbv\ is nonnegative.  

Let \mbx\ be a nonconstant vector in $\R^n_{>0}$.  
Applying Lemma~\ref{Muirhead:lemma:transposition}, we have 
\begin{align*}
[\mbx^{\mbu}]_G 
 & = \frac{1}{2|G|} \sum_{\sigma \in G} 
 \left( x_{\sigma  (j)}^{u_j} x_{\sigma (k)}^{u_{k} }  +  x_{\sigma  (j)}^{u_{k}}  x_{\sigma (k)}^{u_{j} }  \right) 
\prod_{\substack{i=1 \\ i \neq j,k }}^n  x_{\sigma(i)}^{u_i} \\ 
& = \frac{1}{2|G|} \sum_{\sigma \in G} 
 \left( x_{\sigma  (j)}^{\rho + \Delta} x_{\sigma (k)}^{\rho - \Delta}  
 +  x_{\sigma  (j)}^{\rho - \Delta}  x_{\sigma (k)}^{\rho + \Delta}  \right) 
\prod_{\substack{i=1 \\ i \neq j,k }}^n  x_{\sigma(i)}^{u_i}  
\end{align*}
and 
\begin{align*}
[\mbx^{\mbv}]_G 
& = \frac{1}{2|G|} \sum_{\sigma \in G} 
 \left( x_{\sigma  (j)}^{v_j} x_{\sigma (k)}^{v_k}  
 +  x_{\sigma  (j)}^{v_k}  x_{\sigma (k)}^{v_j}  \right) 
\prod_{\substack{i=1 \\ i \neq j,k }}^n  x_{\sigma(i)}^{v_i} \\ 
& = \frac{1}{2|G|} \sum_{\sigma \in G} 
 \left( x_{\sigma  (j)}^{\rho + \delta} x_{\sigma (k)}^{\rho - \delta}  
 +  x_{\sigma  (j)}^{\rho - \delta}  x_{\sigma (k)}^{\rho + \delta}  \right) 
\prod_{\substack{i=1 \\ i \neq j,k }}^n  x_{\sigma(i)}^{u_i}. 
\end{align*}
Therefore,
\begin{align*}
& [\mbx^{\mbu}]_G - [\mbx^{\mbv}]_G \\
& = \frac{1}{2|G|} \sum_{\sigma \in G} 
 \left( 
 x_{\sigma  (j)}^{\rho + \Delta} x_{\sigma (k)}^{\rho - \Delta}  
 +  x_{\sigma  (j)}^{\rho - \Delta}  x_{\sigma (k)}^{\rho + \Delta} 
  - x_{\sigma  (j)}^{\rho + \delta} x_{\sigma (k)}^{\rho - \delta}  
 -  x_{\sigma  (j)}^{\rho - \delta}  x_{\sigma (k)}^{\rho + \delta}  
 \right) 
\prod_{\substack{i=1 \\ i \neq j,k }}^n  x_{\sigma(i)}^{u_i}. 
\end{align*}
If $x_{\sigma(j)} = x_{\sigma(k)}$, then 
\[
 x_{\sigma  (j)}^{\rho + \Delta} x_{\sigma (k)}^{\rho - \Delta}  
 +  x_{\sigma  (j)}^{\rho - \Delta}  x_{\sigma (k)}^{\rho + \Delta} 
  - x_{\sigma  (j)}^{\rho + \delta} x_{\sigma (k)}^{\rho - \delta}  
 -  x_{\sigma  (j)}^{\rho - \delta}  x_{\sigma (k)}^{\rho + \delta}   = 0. 
\]
By~\eqref{Muirhead:SimpleIneq-3} of Corollary~\ref{Muirhead:corollary:SimpleIneq-2}, 
if  $x_{\sigma(j)} \neq x_{\sigma(k)}$, then 
\[
 x_{\sigma  (j)}^{\rho + \Delta} x_{\sigma (k)}^{\rho - \Delta}  
 +  x_{\sigma  (j)}^{\rho - \Delta}  x_{\sigma (k)}^{\rho + \Delta} 
  - x_{\sigma  (j)}^{\rho + \delta} x_{\sigma (k)}^{\rho - \delta}  
 -  x_{\sigma  (j)}^{\rho - \delta}  x_{\sigma (k)}^{\rho + \delta} > 0. 
\]
Because $x_{\sigma(j)} \neq x_{\sigma(k)}$ for some $\sigma \in G$, 
we have  $[\mbx^{\mbu}]_G - [\mbx^{\mbv}]_G > 0$.  
This completes the proof. 
\end{proof}

\section{Muirhead's inequality}

The vector  $\mba = \vectorsmallan \in \Rn$ 
is \emph{decreasing}\index{vector!decreasing} if 
\[
a_1 \geq a_2 \geq  \cdots \geq a_n.
\]
Associated to every vector $\mba = \vectorsmallan \in \Rn$ is a unique decreasing vector 
$\mba^{\downarrow} = \bsmallmat a_1^{\downarrow}  \\ \vdots \\ a_n^{\downarrow}\esmallmat \in \Rn $ 
obtained from \mba\ by a rearrangement of  coordinates.  
Thus, $\mba^{\downarrow} = \sigma (\mba)$ for some $\sigma \in S_n$. 

For example,  from the vector 
$\mba =  \bsmallmat 3 \\ 1 \\ 2 \\ 4 \esmallmat \in \R^4$ we obtain the decreasing vector
$\mba^{\downarrow} = \bsmallmat 4 \\ 3 \\ 2 \\ 1 \esmallmat$.  The permutation 
\[
\sigma = \bmat 1 & 2 & 3 & 4 \\ 2 & 4 & 3 & 1 \emat \in S_4  
\]
satisfies $\mba^{\downarrow} = \sigma(\mba)$.  

Let $\mba = \vectorsmallan$ and $\mbb = \vectorsmallbn$ be vectors in \Rn, 
and let $\mba^{\downarrow} = \bsmallmat a_1^{\downarrow} \\ \vdots \\ a_n^{\downarrow} \esmallmat$ 
and $\mbb^{\downarrow} = \bsmallmat b_1^{\downarrow} \\ \vdots \\ b_n^{\downarrow} \esmallmat$ 
be the corresponding decreasing vectors obtained by permutation of coordinates.    
The following definition is fundamental: 
The vector \mba\ \emph{majorizes} the vector \mbb, denoted $\mbb \prec \mba$, if 
\[
\sum_{i=1}^k b_i \leq \sum_{i=1}^k a_i \qquad \text{for all $i \in \{1,\ldots, n-1\}$ }
\]
and
\[
\sum_{i=1}^n b_i = \sum_{i=1}^n a_i.
\]
The vector \mba\ \emph{strictly majorizes} \mbb\ if $\mbb \prec \mba$  and $\mbb \neq \mba$.

The \emph{Hamming distance} between vectors $\mba = \vectorsmallan$ 
and $\mbb =\vectorsmallbn$ in $\R^n$ is
\beq        \label{Muirhead:Hamming}
d_H(\mba,\mbb) = \card\{i \in \{ 1,\ldots, n\}: a_i  \neq b_i \}.
\eeq

\bl              \label{Muirhead:lemma:Muirhead-transposition-2}
Let $\mbu = \vectorsmallun$ and $\mbv = \vectorsmallvn$ be vectors in $\Rn_{\geq 0}$  
such that $\mbv \prec \mbu$ and $d_H(\mbu,\mbv) = 2$.  Then 
\[
\left[ \mbx^{\mbv} \right]_{S_n} < \left[ \mbx^{\mbu} \right]_{S_n} 
\]
for all nonconstant vectors $\mbx \in \R^n_{>0}$.   
\el

\begin{proof}
By Lemma~\ref{Muirhead:lemma:a-permute}, we can assume that the vectors \mbu\ 
and \mbv\ are decreasing.  

Hamming distance $d_H(\mbu,\mbv) = 2$ implies that there are unique integers $j < k$ 
such that $u_i = v_i$ for all $i \neq j,k$ and 
\[
0 \leq u_k < v_k  \leq v_j < u_j \qqand v_j + v_k = u_j + u_k.
\]
The numbers 
\[
\rho = \frac{u_j + u_k}{2} = \frac{v_j + v_k}{2}, \qquad \Delta = \frac{u_j - u_k}{2} 
\]
and 
\[
\delta = v_j - \frac{u_j + u_k}{2} = \frac{v_j - v_k}{2}
\]
satisfy the inequality 
\[
0 \leq \delta < \Delta \leq \rho.  
\] 
We have 
\[
 \rho + \delta =  \frac{u_j + u_k}{2} + \left( v_j - \frac{u_j + u_k}{2}   \right) =  v_j 
\]
and 
\[
\rho - \delta = \frac{v_j + v_k}{2} - \frac{v_j - v_k}{2} =  v_k.  
\]
Applying Lemma~\ref{Muirhead:lemma:Muirhead-transposition} 
completes the proof. 
\end{proof}

\bt                                     \label{Muirhead:theorem:HardyLittlewoodPolya-interpolation}
Let \mba\ and \mbb\ be distinct  vectors in $\Rn$ such that $\mbb \prec \mba$.
For some positive integer  $r \leq d_H(\mba,\mbb)$, 
there is a sequence of decreasing vectors 
$\mbc_0, \mbc_1, \ldots, \mbc_r \in \R^n_{\geq 0}$ such that 
\beq                   \label{Muirhead:HardyLittlewoodPolya-interpolation-1}
\mbb = \mbc_r \prec \mbc_{r-1} \prec \cdots \prec \mbc_1 \prec \mbc_0 = \mba 
\eeq
and 
\beq                   \label{Muirhead:HardyLittlewoodPolya-interpolation-2}
d_H(\mbc_{i-1}, \mbc_i) = 2 \qquad \text{for all $i \in \{1,\ldots, r\}$.}
\eeq
\et

\begin{proof}
This is Theorem 4 in  Nathanson~\cite{nath22}.
\end{proof}

We observe that if $\mba = \vectorsmallan$ and $\mbb = \vectorsmallbn$ 
and if $\sum_{i=1}^n a_i = \sum_{i=1}^n b_i$, then,  
for every $x > 0$, the constant vector 
$\mbx  = \bsmallmat x \\ \vdots \\ x \esmallmat \in \R^n_{>0}$ 
satisfies  
\[
[\mbx^{\mbb}]_{S_n} = [\mbx^{\mba}]_{S_n} = x^{\sum_{i=1}^n a_i}.
\] 
In particular, if $\mbb \prec \mba$, then $[\mbx^{\mbb}]_{S_n} = [\mbx^{\mba}]_{S_n}$ 
for every constant vector \mbx.
Muirhead's inequality is a strict inequality for nonconstant vectors.

\bt[Muirhead's inequality]         \label{Muirhead:theorem:Muirhead}
Let \mba\ and \mbb\ be distinct vectors in $\R^n_{\geq 0}$.   
If $\mbb \prec \mba$, then 
\[
[\mbx^{\mbb}]_{S_n} < [\mbx^{\mba}]_{S_n}
\]
for every nonconstant vector $\mbx\in \R^n_{>0}$. 
\et

\begin{proof} 
We give two proofs of this inequality.  
By Theorem~\ref{Muirhead:theorem:HardyLittlewoodPolya-interpolation}, 
 there are decreasing vectors $\mbc_1,\ldots, \mbc_r$ 
 and a strict majorization chain  
\[
\mbb = \mbc_r \prec \mbc_{r-1} \prec \cdots \prec \mbc_1 \prec \mbc_0 = \mba 
\]
such that 
\[
d_H(\mbc_{i-1}, \mbc_i) = 2 \qquad \text{for all $i \in \{1,\ldots, r\}$.}
\]

Let \mbx\ be a nonconstant vector in $\R^n_{>0}$.  
The symmetric group $S_n$ contains all transpositions, and so,   
by Lemma~\ref{Muirhead:lemma:Muirhead-transposition-2} 
with $\mbu = \mbc_{i-1}$ and $\mbv = c_i$,  
\[
[\mbx^{\mbc_i}]_{S_n} < [\mbx^{\mbc_{i-1}}]_{S_n}
\]
for all $i \in \{1,\ldots, r\}$.  Therefore,  
\[
[\mbx^{\mbb}]_{S_n} = [\mbx^{\mbc_r}]_{S_n} <   [\mbx^{\mbc_{r-1}}]_{S_n} 
< \cdots <  [\mbx^{\mbc_1}]_{S_n}<  [\mbx^{\mbc_0}]_{S_n} = [\mbx^{\mba}]_{S_n}. 
\]
This completes the first proof.  
\end{proof} 

The second proof of Theorem~\ref{Muirhead:theorem:Muirhead} 
uses only the arithmetic and  geometric mean inequality.

\begin{proof} 
Let $(c_i)_{i=1}^n$ be a sequence of real numbers such that $c_i \neq 0$ for some $i$ and 
\[
\sum_{i=1}^n c_i = 0.
\]
Let $\mbx = \vectorsmallxn$ be a nonconstant positive vector.  
For every pair of integers $(i,j)$ with $i,j \in \{1,\ldots, n\}$, there are $(n-1)!$ permutations 
$\sigma \in S_n$ with $\sigma(i) = j$ and so there are $(n-1)!$ permutations $\sigma$ 
such that $x_{\sigma(i)}^{c_i} = x_j^{c_i}$.  
Therefore,
\begin{align*}
\prod_{\sigma \in S_n}  x_{\sigma(1)}^{c_1} x_{\sigma(2)}^{c_2} \cdots x_{\sigma(n)}^{c_n} 
& = \prod_{\sigma \in S_n} \prod_{i=1}^n  x_{\sigma(i)}^{c_i} 
= \prod_{i=1}^n \prod_{\sigma \in S_n}  x_{\sigma(i)}^{c_i}  \\
& = \prod_{i=1}^n \prod_{j =1}^n  x_j^{(n-1)! c_i} 
= \prod_{j =1}^n \prod_{i =1}^n x_j^{(n-1)! c_i} \\
& =  \prod_{j =1}^n   x_j^{(n-1)! \sum_{i=1}^n c_i} = 1.
\end{align*}
The arithmetic and geometric mean inequality  gives 
\[
1 = \left( \prod_{\sigma \in S_n} x_{\sigma(1)}^{c_1} x_{\sigma(2)}^{c_2} \cdots x_{\sigma(n)}^{c_n} \right)^{1/n!} 
< \frac{1}{n!}\sum_{\sigma \in S_n}  x_{\sigma(1)}^{c_1} x_{\sigma(2)}^{c_2} \cdots x_{\sigma(n)}^{c_n}
\]
and so 
\[
n! < \sum_{\sigma \in S_n}  x_{\sigma(1)}^{c_1} x_{\sigma(2)}^{c_2} \cdots x_{\sigma(n)}^{c_n}.
\]
Equivalently, 
 \beq               \label{Muirhead:constantMuirhead}
0 < \sum_{\sigma \in S_n}  \left( x_{\sigma(1)}^{c_1} x_{\sigma(2)}^{c_2} \cdots x_{\sigma(n)}^{c_n} - 1\right).
\eeq

The vector $\mba = \vectorsmallan$ strictly majorizes the vector $\mbb = \vectorsmallbn$ 
and so $\mba \neq \mbb$.  Let
\[
c_i = a_i-b_i 
\]
 for $i \in\{1,\ldots, n\}$.  Then 
\[
\sum_{i=1}^n c_i = \sum_{i=1}^n (a_i-b_i) = \sum_{i=1}^n a_i - \sum_{i=1}^n b_i  = 0 
\]
and $c_i \neq 0$ for some $i \in\{1,\ldots, n\}$.    
We have 
\begin{align*}
\sum_{\sigma \in S_n} \prod_{i=1}^n x_{\sigma(i)}^{a_i} 
& = \sum_{\sigma \in S_n} \prod_{i=1}^n x_{\sigma(i)}^{b_i+c_i} 
 = \sum_{\sigma \in S_n} \prod_{i=1}^n x_{\sigma(i)}^{b_i} \prod_{i=1}^n x_{\sigma(i)}^{c_i}. 
\end{align*} 
Inequality~\eqref{Muirhead:constantMuirhead} implies   
\begin{align*}
[\mbx^{\mba}]_{S_n} - [\mbx^{\mbb}]_{S_n}
& = \frac{1}{n!}  \sum_{\sigma \in S_n} \prod_{i=1}^n x_{\sigma(i)}^{a_i} 
 - \frac{1}{n!}  \sum_{\sigma \in S_n} \prod_{i=1}^n x_{\sigma(i)}^{b_i} \\
& = \frac{1}{n!} \sum_{\sigma \in S_n} \prod_{i=1}^n x_{\sigma(i)}^{b_i} \prod_{i=1}^n x_{\sigma(i)}^{c_i} 
 - \frac{1}{n!} \sum_{\sigma \in S_n} \prod_{i=1}^n x_{\sigma(i)}^{b_i} \\
 & = \frac{1}{n!} \sum_{\sigma \in S_n} \prod_{i=1}^n x_{\sigma(i)}^{b_i} 
 \left( \prod_{i=1}^n x_{\sigma(i)}^{c_i} -1\right)  \\
 & > 0.  
\end{align*} 
This completes the proof.  
\end{proof}

Next we prove the converse of Muirhead's inequality.

\bt         \label{Muirhead:theorem:Muirhead-2}
Let $\mba = \vectorsmallan$ and $\mbb = \vectorsmallbn$ 
be  vectors in $\R^n_{\geq 0}$. 
If 
\[
[\mbx^{\mbb}]_{S_n} \leq [\mbx^{\mba}]_{S_n}
\]
for every  vector $\mbx\in \R^n_{>0}$, 
then $\mbb \prec \mba$.   
\et

\begin{proof}  
For  all $x > 0$, the constant vector $\mbx = \bsmallmat x \\ \vdots \\ x \esmallmat$ 
is positive and satisfies   
\[
\left[\mbx^{\mba} \right]_{S_n} =  x^{\sum_{i=1}^n a_i} 
\qqand 
\left[\mbx^{\mbb} \right]_{S_n} 
=  x^{\sum_{i=1}^n b_i}.
\]
Suppose that  $[\mbx^{\mbb}]_{S_n}  \leq [\mbx^{\mba}]_{S_n}$ for every 
positive vector $\mbx$.  For all $x > 1$ we have 
\[
x^{\sum_{i=1}^n b_i} \leq x^{\sum_{i=1}^n a_i} 
\]
and so 
\[
\sum_{i=1}^n b_i \leq \sum_{i=1}^n a_i.
\]
Similarly, $1/x > 0$ and 
\[
\frac{1}{x^{\sum_{i=1}^n b_i}} = \left(\frac{1}{x}\right)^{\sum_{i=1}^n b_i} 
\leq \left(\frac{1}{x}\right)^{\sum_{i=1}^n a_i} 
= \frac{1}{x^{\sum_{i=1}^n a_i}} 
\]
and so 
\[
\sum_{i=1}^n b_i \geq \sum_{i=1}^n a_i.
\]
Therefore, 
\[
\sum_{i=1}^n b_i = \sum_{i=1}^n a_i.
\]

By Lemma~\ref{Muirhead:lemma:a-permute}, 
we can assume that the vectors \mba\ and \mbb\ are decreasing.  
For  $k \in \{1,\ldots, n-1\}$, we define  the variables 
\[
x_i = \begin{cases}
x & \text{if $i \in \{1,\ldots, k\}$ }\\
1 & \text{if $i \in \{k+1,\ldots, n\}$} 
\end{cases}
\] 
 the vector 
\[
\mbx = \bsmallmat x_1 \\ \vdots \\ x_k \\ x_{k+1} \\ \vdots \\ x_n \esmallmat 
= \bsmallmat x  \\ \vdots \\ x  \\1\\ \vdots \\ 1 \esmallmat 
\]
and  the polynomial 
\[
\left[\mbx^{\mba}\right]_{S_n} = \frac{1}{n!} \sum_{\sigma \in S_n} x_{\sigma(1)}^{a_1} \cdots x_{\sigma(n)}^{a_n}. 
\]
The leading coefficient of the polynomial $\left[\mbx^{\mba}\right]_{S_n} $ is positive 
and, because the vector \mba\ is decreasing,  the degree of this polynomial is 
$ \sum_{i=1}^k a_i$.  
Similarly, the polynomial 
\[
\left[\mbx^{\mbb}\right]_{S_n}  
= \frac{1}{n!} \sum_{\sigma \in S_n} x_{\sigma(1)}^{b_1} \cdots x_{\sigma(n)}^{b_n}
\]
has a positive leading coefficient and degree $ \sum_{i=1}^k b_i$.
The inequality $\left[\mbx^{\mbb}\right]_{S_n} \leq \left[\mbx^{\mba}\right]_{S_n}$ 
for all $x > 0$ 
implies that $\sum_{i=1}^k b_i \leq \sum_{i=1}^k a_i$ for all $k \in \{1,2,\ldots, n\}$  
and so $\mbb \prec \mba$.  
This completes the proof. 
\end{proof}

We use Muirhead's inequality to prove the following result.

\bt                                  \label{Muirhead:theorem:Muirhead-AddMult}
Let $(u_i)_{i=1}^n$ and $(v_i)_{i=1}^n$ be decreasing sequences of positive numbers.     
If $(u_i)_{i=1}^n \neq (v_i)_{i=1}^n$ and 
\beq                                    \label{Muirhead:Muirhead-AddMult-Ineq1}
\prod_{i=1}^j v_i \leq \prod_{i=1}^j u_i \qquad\text{for all $j \in \{1,\ldots, n\}$}
\eeq
then 
\beq                                \label{Muirhead:Muirhead-AddMult-Ineq2}
\sum_{i=1}^n v_i < \sum_{i=1}^n u_i.
\eeq
\et

\begin{proof} 
Let 
\[
v_{n+1} = \min(v_n,u_n)  \qqand 
u_{n+1} = v_{n+1}  \left( \frac{\prod_{i=1}^n v_i}{\prod_{i=1}^n u_i}\right). 
\]
We have $ v_n \geq v_{n+1} > 0$ and so the sequence  $(v_i)_{i=1}^{n+1}$ 
is positive and decreasing.   
Inequality~\eqref{Muirhead:Muirhead-AddMult-Ineq1} with $j=n$ implies that 
$\prod_{i=1}^n u_i \geq \prod_{i=1}^n v_i > 0$ and  so 
\[
0 < u_{n+1} = v_{n+1}   \left( \frac{\prod_{i=1}^n v_i}{\prod_{i=1}^n u_i} \right) 
\leq v_{n+1} =  \min(v_n,u_n) \leq u_n.
\]
Thus, the sequence  $(u_i)_{i=1}^{n+1}$ is also positive and decreasing. 
Moreover, 
\begin{align*}
\prod_{i=1}^{n+1} u_i 
& = u_{n+1}\prod_{i=1}^n u_i 
= v_{n+1} \left( \frac{\prod_{i=1}^n v_i}{\prod_{i=1}^n u_i}\right) \prod_{i=1}^n u_i 
 = \prod_{i=1}^{n+1} v_i.
\end{align*}

Let $\lambda > 0$.  For all  $j \in \{1,\ldots, n\}$, we have  
\[
\prod_{i=1}^j \lambda v_i \leq \prod_{i=1}^j \lambda u_i \qquad\text{if and only if} \qquad \prod_{i=1}^j v_i \leq \prod_{i=1}^j u_i 
\]
and
\[
\sum_{i=1}^j \lambda v_i \leq \sum_{i=1}^j \lambda u_i \qquad\text{if and only if} \qquad \sum_{i=1}^j v_i \leq \sum_{i=1}^j u_i. 
\]
Choosing $\lambda$ sufficiently large, we can assume that $v_i > 1$ 
and $u_i >  1$ for all $i \in \{1,\ldots, n, n+1\}$, and so 
\[
 a_i = \log u_i  >  0  \qqand b_i = \log v_i  > 0 
\]
for all $i \in \{1,\ldots, n, n+1\}$.  
The sequences $(b_i)_{i=1}^{n+1}$ and $(a_i)_{i=1}^{n+1}$ 
are decreasing sequences of positive numbers.  
Moreover, $(u_i)_{i=1}^{n+1} \neq (v_i)_{i=1}^{n+1}$ implies 
$(a_i)_{i=1}^{n+1} \neq (b_i)_{i=1}^{n+1}$. 

Let 
\[
\mba =  \bmat a_1 \\ \vdots \\ a_n \\ a_{n+1} \emat 
\qqand 
\mbb = \bmat b_1 \\ \vdots \\ b_n \\ b_{n+1} \emat.  
\]
For all $j \in \{1,\ldots, n\}$ we have 
\[
\sum_{i=1}^j b_i = \sum_{i=1}^j \log v_i = \log\prod_{i=1}^j v_i 
\leq \log\prod_{i=1}^j u_i = \sum_{i=1}^j \log u_i = \sum_{i=1}^j a_i 
\]
and 
\[
\sum_{i=1}^{n+1} b_i = \sum_{i=1}^{n+1} \log v_i = \log\prod_{i=1}^{n+1} v_i 
= \log\prod_{i=1}^{n+1} u_i =  \sum_{i=1}^{n+1} \log u_i = \sum_{i=1}^{n+1} a_i.
\]
Thus, the vector \mba\ strictly majorizes the vector \mbb. 

By Muirhead's inequality (Theorem~\ref{Muirhead:theorem:Muirhead}), for every 
positive vector $\mbx = \vectorsmallxn$, we have 
\begin{align*}
\frac{1}{n!} \sum_{\sigma \in S_{n+1}} \prod_{i=1}^{n+1} x_{\sigma(i)}^{\log v_i}  
& = \frac{1}{n!} \sum_{\sigma \in S_{n+1}}  \prod_{i=1}^{n+1} x_{\sigma(i)}^{b_i} 
  = \left[ \mbx^{\mbb}\right]_{S_n} \\ 
& <   \left[ \mbx^{\mba}\right]_{S_n}          
 = \frac{1}{n!} \sum_{\sigma \in S_{n+1}}  \prod_{i=1}^{n+1} x_{\sigma(i)} ^{a_i} \\
&  =  \frac{1}{n!}  \sum_{\sigma \in S_{n+1}} \prod_{i=1}^{n+1} x_{\sigma(i)} ^{\log u_i}.  
\end{align*}
Let  $x_1 = e$  and $x_i = 1$ for $i \in \{2,\ldots, n+1\}$. 
If $\sigma \in S_{n+1}$ and $\sigma(j) = 1$, then 
\[
  \prod_{i=1}^{n+1} x_{\sigma(i)}^{b_i}  =  \prod_{i=1}^{n+1} x_{\sigma(i)}^{\log v_i} 
  = x_1^{\log v_j}   = e^{\log v_j} = v_j.
 \]
There are $n!$ permutations $\sigma \in S_{n+1}$ such that $\sigma(j) = 1$, 
and so  
\[
n! v_j =  \sum_{\substack{\sigma \in S_{n+1} \\ \sigma(j) = 1}} v_j. 
\]
Similarly, 
\[
n! u_j =  \sum_{\substack{\sigma \in S_{n+1} \\ \sigma(j) = 1}} u_j. 
\]
Therefore, 
\begin{align*}
n! \sum_{j=1}^{n+1} v_j  
& = \sum_{j=1}^{n+1} \sum_{\substack{\sigma \in S_{n+1} \\ \sigma(j) = 1}} v_j   
 = \sum_{j=1}^{n+1} \sum_{\substack{\sigma \in S_{n+1} \\ \sigma(j) = 1}} x_1^{b_j}   \\
& =  \sum_{j=1}^{n+1} \sum_{\substack{\sigma \in S_{n+1} \\ \sigma(j) = 1}} 
 \prod_{i=1}^{n+1} x_{\sigma(i)}^{b_i} 
  =   \sum_{\sigma \in S_{n+1}}  \prod_{i=1}^{n+1} x_{\sigma(i)}^{b_i}  \\
& <  \sum_{\sigma \in S_{n+1}}  \prod_{i=1}^{n+1} x_{\sigma(i)}^{a_i} 
= \sum_{i=1}^{n+1} \sum_{\substack{\sigma \in S_{n+1} \\ \sigma(j) = 1}}  
 \prod_{i=1}^{n+1} x_{\sigma(i)} ^{a_i} \\
 & = \sum_{j=1}^{n+1} \sum_{\substack{\sigma \in S_{n+1} \\ \sigma(j) = 1}} x_1^{a_j}   = \sum_{j=1}^{n+1} \sum_{\substack{\sigma \in S_{n+1} \\ \sigma(j) = 1}} u_j   \\ 
 & = n! \sum_{j=1}^{n+1} u_j.
\end{align*}
The inequality $u_{n+1} \leq  v_{n+1}$ implies  
\[
\sum_{i=1}^{n} v_i  = \sum_{i=1}^{n+1} v_i - v_{n+1} 
<  \sum_{i=1}^{n+1} u_i - u_{n+1} = \sum_{i=1}^n u_i.  
\] 
This completes the proof.  
\end{proof}

\section{Rado's inequality} 

Let $G$ be a subgroup of $S_n$ and let  $\mba$ be 
a vector in $\R^n$.  
The \emph{$G$-permutohedron}\index{permutohedron} generated by  \mba,  
denoted $K_{G}(\mba)$,  
is the convex hull of the finite set of vectors 
$\{\gamma(\mba):\gamma \in G\}$. 
The $G$-permutohedron $K_G(\mba)$ is a  compact convex subset of \Rn.  
If the vector $\mba$ is nonnegative, then $\gamma(\mba)$ 
is nonnegative for all $\gamma \in G$ 
and so every convex combination of the vectors $\gamma(\mba)$, 
that is, every vector   in $K_G(\mba)$, is nonnegative. 

We recall the following separation theorem from convexity theory.

\bt                                     \label{Muirhead:theorem:hyperplane}
Let $K$ be a compact convex set 
in \Rn.  If the vector $\mbb \in \Rn$ is not in $K$, then the set $K$ 
and the vector $\mbb$ are 
strictly separated by a hyperplane.
\et

This means that there is a nonzero linear functional
\[
H(\mbx) = \sum_{i=1}^n u_i x_i 
\]
and there are scalars $c$ and  $\delta$  with $\delta > 0$ such that 
\[
H(\mbx) \leq c \qquad \text{for all $\mbx \in K$}
\] 
and 
\[
H(\mbb) \geq c + \delta. 
\]

We also use the general form of the \emph{arithmetic-geometric mean inequality}.

\bt                     \label{Muirhead:theorem:AGM}
Let $G$ be a finite set.  If $\{w_{\gamma}:\gamma \in G\}$ and 
$\{t_{\gamma}:\gamma \in G \}$ are sets of nonnegative numbers such that 
$\sum_{ \gamma \in G} t_{\gamma} = 1$, then 
\[
\prod_{\gamma \in G} w_{\gamma}^{t_{\gamma}} \leq \sum_{\gamma \in G} t_{\gamma} w_{\gamma}. 
\]
Moreover, 
\[
\prod_{\gamma \in G} w_{\gamma}^{t_{\gamma}} = \sum_{\gamma \in G} t_{\gamma} w_{\gamma}. 
\]
if and only if $w_{\gamma} = w_{\gamma'}$ for all $\gamma,\gamma' \in G$.  
\et
The following result of Richard Rado extends Muirhead's inequality to subgroups of $S_n$.

\bt[Rado]                 \label{Muirhead:theorem:Rado-1}
Let $G$ be a subgroup of the symmetric group $S_n$.  
Let $\mba$  be a nonnegative vector in $\Rn$ and let $K_G(\mba)$ be the $G$-permutohedron, 
that is, the convex hull of the set of vectors $\{\gamma (\mba): \gamma \in G\}$.  
If 
\[
 \mbb \in K_G(\mba)
\]
then $\mbb$ is nonnegative and 
\[
\left[\mbx^{\mbb}\right]_G \leq \left[\mbx^{\mba}\right]_G 
\]
for all positive vectors $\mbx  \in \Rn$. 
\et

\begin{proof} 
Let $\mba = \vectorsmallan \in \R^n_{\geq 0}$ and $\mbb = \vectorsmallbn \in K_G(\mba)$. 
Then 
\[
\mbb = \sum_{\gamma \in G} t_{\gamma} \ \gamma (\mba)
\]
where $t_{\gamma} \in [0,1]$ for all $\gamma \in G$ and 
$\sum_{\gamma \in G} t_{\gamma} = 1$.  
By~\eqref{Muirhead: sigma-z},  the $i$th component of the vector \mbb\ is  
\[
b_i = \sum_{\gamma \in G} t_{\gamma} (\gamma \mba)_i  = \sum_{\gamma \in G} t_{\gamma} \ a_{\gamma ^{-1}(i)}. 
\]
Let $\mbx = \vectorsmallxn$ be a positive vector.   
For every permutation $\sigma \in G$, the arithmetic  and geometric mean inequality 
(Theorem~\ref{Muirhead:theorem:AGM}) gives 
\[
 \prod_{\gamma \in G} \left( \prod_{i=1}^n x_{\sigma(i)}^{ a_{\gamma ^{-1}(i)} } \right)^{ t_{\gamma}} \\
\leq   \sum_{\gamma \in G}  t_{\gamma} \left( \prod_{i=1}^n x_{\sigma(i)}^{ a_{\gamma ^{-1}(i)} } \right)  
\]
and so 
\begin{align*}
\left[\mbx^{\mbb}\right]_G 
& = \frac{1}{|G|} \sum_{\sigma \in G} \prod_{i=1}^n x_{\sigma(i)}^{b_i} 
 = \frac{1}{|G|} \sum_{\sigma \in G} \prod_{i=1}^n x_{\sigma(i)}^{ \sum_{\gamma \in G} t_{\gamma} a_{\gamma ^{-1}(i)} } \\ 
& = \frac{1}{|G|} \sum_{\sigma \in G} \prod_{i=1}^n  \prod_{\gamma \in G} x_{\sigma(i)}^{ t_{\gamma} a_{\gamma ^{-1}(i)} } 
  = \frac{1}{|G|} \sum_{\sigma \in G}  \prod_{\gamma \in G} \left( \prod_{i=1}^n x_{\sigma(i)}^{ a_{\gamma ^{-1}(i)} } \right)^{ t_{\gamma}} \\
& \leq   \frac{1}{|G|} \sum_{\sigma \in G}  \sum_{\gamma \in G}  t_{\gamma} \left( \prod_{i=1}^n x_{\sigma(i)}^{ a_{\gamma ^{-1}(i)} } \right)   
 = \frac{1}{|G|} \sum_{\gamma \in G}  t_{\gamma}   \sum_{\sigma \in G}  \left( \prod_{i=1}^n x_{\sigma(i)}^{ a_{\gamma ^{-1}(i)} } \right) \\
& = \frac{1}{|G|} \sum_{\gamma \in G}  t_{\gamma}   \sum_{\sigma \in G}  \left( \prod_{j=1}^n x_{\sigma \gamma (j)}^{ a_j } \right) 
= \frac{1}{|G|} \sum_{\gamma \in G}  t_{\gamma}   \sum_{\sigma \in G}  \left( \prod_{j=1}^n x_{\sigma  (j)}^{ a_j } \right)  \\
&  =  \frac{1}{|G|} \sum_{\sigma \in G}  \left( \prod_{j=1}^n x_{\sigma  (j)}^{ a_j } \right) \sum_{\gamma \in G}  t_{\gamma}   
=  \frac{1}{|G|} \sum_{\sigma \in G}  \left( \prod_{j=1}^n x_{\sigma  (j)}^{ a_j } \right) \\ 
& = \left[\mbx^{\mba}\right]_G.
\end{align*} 
This completes the proof.  
\end{proof}

Next we prove the converse of Theorem~\ref{Muirhead:theorem:Rado-1}. 

\bt[Rado]                        \label{Muirhead:theorem:Rado-2}
Let $G$ be a subgroup of the symmetric group $S_n$.  
Let $\mba$  be a nonnegative vector in $\Rn$ and let $K_G(\mba)$ be the $G$-permutohedron. 
If  \mbb\ is a nonnegative vector in $\Rn$ such that 
\[
\left[\mbx^{\mbb}\right]_G \leq \left[\mbx^{\mba}\right]_G 
\]
for all positive vectors $\mbx  \in \Rn$, then 
\[
 \mbb \in K_G(\mba). 
\]
\et

\begin{proof}
Let $\mba = \vectorsmallan$ and $\mbb = \vectorsmallbn$.  
We shall prove that if $\mbb \notin K_G(\mba)$, then there is a positive 
vector $\mbx \in \Rn$ such that $\left[\mbx^{\mbb}\right]_G > \left[\mbx^{\mba}\right]_G$. 

By Theorem~\ref{Muirhead:theorem:hyperplane}, 
if $\mbb \notin K_G(\mba)$, then the compact convex set $K_G(\mba)$ 
and the vector $\mbb$ are 
strictly separated by a hyperplane $H$ and so there is a nonzero linear functional
\[
H(\mbx) = \sum_{i=1}^n u_i x_i 
\]
and scalars $c$ and  $\delta$  with $\delta > 0$ such that 
\[
H(\mbx) \leq c \qquad \text{for all $\mbx \in K_G(\mba)$}
\] 
and 
\[
H(\mbb) \geq c + \delta. 
\]
For all $\gamma \in G$ we have 
\[
\gamma^{-1} (\mba)  = \bmat a_{\gamma(1)} \\ \vdots \\ a_{\gamma(n)} \emat \in K_G(\mba)
\]
 and so 
\[
 \sum_{i=1}^n u_i a_{\gamma^{-1}(i)} = H(\gamma (\mba))   \leq  c  
 \leq  H(\mbb) -  \delta = \sum_{i=1}^n u_i b_i -  \delta. 
\]
  
Let 
\[
M > |G|^{1/\delta} > 1 \qqand x_i = M^{u_i} 
\]
for all $i \in \{1,\ldots, n\}$.  
 The vector $\mbx = \vectorsmallxn$  is positive. 
 The subgroup $G$ contains the identity permutation, and so 
\begin{align*} 
M^{\sum_{i=1}^n  u_i b_i} 
& \leq  \sum_{\gamma \in G} M^{\sum_{i=1}^n  u_i b_{\gamma(i)}} = \sum_{\gamma \in G} \prod_{i=1}^n M^{u_i b_{\gamma(i)}} \\
 & = \sum_{\gamma \in G} \prod_{i=1}^n x_i^{b_{\gamma(i)}} 
= |G| \left[\mbx^\mbb \right]_G. 
\end{align*} 
We have $\gamma(\mba) \in K_G(\mba)$ and $H(\gamma(\mba)) \geq H(\mbb)-\delta$.  
It follows that 
\begin{align*}
 \left[\mbx^\mba\right]_G
& = \frac{1}{|G|}  \sum_{\gamma \in G} \prod_{i=1}^n x_i^{a_{\gamma^{-1}(i)}} 
=   \frac{1}{|G|} \sum_{\gamma \in G} \prod_{i=1}^n M^{u_i a_{\gamma^{-1}(i)}} \\
& =   \frac{1}{|G|} \sum_{\gamma \in G} M^{\sum_{i=1}^n  u_i a_{\gamma^{-1}(i)}}  
 =  \frac{1}{|G|} \sum_{\gamma \in G} M^{   H(\gamma (\mba) )  } \\
& \leq  \frac{1}{|G|}  \sum_{\gamma \in G} M^{   H(\mbb) - \delta  } 
=   \frac{1}{|G|} M^{\sum_{i=1}^n   u_i b_i - \delta } \\
& =   \frac{1}{ M^{\delta} }  \frac{1}{|G|}   M^{\sum_{i=1}^n   u_i b_i }  \leq \frac{1}{ M^{\delta} } \left[\mbx^\mbb \right]_G \\
& <  \left[\mbx^\mbb \right]_G. 
\end{align*}   
This completes the proof.  
\end{proof}

\def\cprime{$'$} \def\cprime{$'$} \def\cprime{$'$}
\providecommand{\bysame}{\leavevmode\hbox to3em{\hrulefill}\thinspace}
\providecommand{\MR}{\relax\ifhmode\unskip\space\fi MR }
\providecommand{\MRhref}[2]{%
  \href{http://www.ams.org/mathscinet-getitem?mr=#1}{#2}
}
\providecommand{\href}[2]{#2}

\end{document}